\documentclass[a4paper,fleqn]{cas-dc}

\usepackage[numbers]{natbib}
\usepackage{graphicx,amsmath,amssymb,mathrsfs,subcaption}

\newdefinition{rem}{Remark}

\begin{document}
\shortauthors{D. Boffi, A. Halim, and G. Priyadarshi}

\title{On the effect of different samplings to the solution of parametric PDE Eigenvalue Problems}
\shorttitle{GPR for parametric eigenproblems}

\author[1,2]{Daniele Boffi}
\cormark[1]
\ead{daniele.boffi@kaust.edu.sa}

\affiliation[1]{organization={Department of Applied Mathematics and Computational Sciences, KAUST},
    country={Saudi Arabia}}
\author[1,3]{Abdul Halim}
\ead{abdul.halim@kaust.edu.sa}
\author[4]{Gopal Priyadarshi}
\ead{gopalpriyadarshi8@gmail.com }
\affiliation[2]{organization={Dipartimento di Matematica ``F. Casorati'', University of Pavia},
    country={Italy}}
\affiliation[3]{organization={Department of Mathematics, H.S. College, Munger University},
    country={India}}
    \affiliation[4]{organization={Department of Mathematics, S.M.D. College, Patliputra University},
    country={India}}

\cortext[cor1]{Corresponding author}

\begin{abstract}
In this article we apply reduced order techniques for the approximation of parametric eigenvalue problems. The effect of the choice of sampling points is investigated. Here we use the standard proper orthogonal decomposition technique to obtain the basis of the reduced space and Galerking orthogonal technique is used to get the reduced problem. We present some numerical results and observe that the use of sparse sampling is a good idea for sampling if the dimension of parameter space is high.
\end{abstract}

\begin{keywords}
Parametric PDE Eigenvalue Problem \sep Sparse Sampling \sep Proper Orthogonal Decomposition \sep Reduced Order Modeling
\end{keywords}

\maketitle

\section{Introduction}
PDE eigenvalue problem is an important class of problems in science and engineering. It is common that PDE eigenvalue problems contain parameters; for example in nuclear reactor physics~\cite{Buchanetal13,auto21,multi_gp} the diffusion of free neutrons is modeled by an eigenvalue problem where the eigenvalue represents the criticality of the reactor and the eigenvector represents the flux of the neutron population. Also, in order to calculate the band gap in a photonic crystals~\cite{photonic1}, one has to solve two parametric eigenvalue problems.
The parameters may represent the properties of the materials, the geometry of the domain, or initial/boundary conditions. The approximation of parametric problems involves the solution of costly eigenvalue problems for several choices of parameters. Hence, we use reduced order modeling for solving the parametric problems and we rely on the offline-online paradigm. In the offline stage, we need to solve the eigenvalue probelm for some parameters. Then we form a reduced system using the solutions at the selected parameters. In the online phase, we obtain the solution using the reduced system for a new parameter.

In order to select the parameters in the offline stage, one can use a greedy algorithm based on some optimization criterion or some sampling technique. If the dimension of the parameter space is small, then we can use uniform sampling but if the dimension is high, then the uniform parameter selection is not of interest. In that case, a viable choice is the use of sparse grid sampling like Latin Hypercube sampling, and Smolyak sparse sampling with Clenshaw Curtis points. In this paper, we investigate the approximation of the eigenvalues of parametric PDEs using different sampling techniques for the selection of the sample parameters in the offline stage of the reduced model.  

In~\cite{WithGopal} we have investigated the choice of snapshots for calculating eigenpairs in ROM and found that if we are interested in approximating the first $k$ eigenvalues then we need to build the snapshot matrix by taking into account all the corresponding eigenfunctions or the sum of them. This is an agreement with the theory developed in~\cite{Fumagallietal16,Horgeretal17} and the results presented in~\cite{GermanRagusa19,Buchanetal13,Pau07a,Pau08}.
Other relevant papers in this context are ~\cite{our_time_cont, Moataz,eccomas,datadriven}. Standard references for the approximation of eigenvalue problems are ~\cite{Boffi10, Boffietal13, Boffietal97}.

The article is organized as follows: in Section~\ref{sec:pb} we discuss the settings of the problem; the details on the Galerkin projection-based reduced order modeling (ROM) is presented in Section~\ref{sec:rom}. Numerical results of two parametric problems are presented in Section~\ref{sec:nr}.
\section{Problem description}\label{sec:pb}
\label{se:ps}
We discuss the effect of the choice of sample points on the solution of the parametric eigenvalue problems using reduced order modeling. Here we consider Galerkin projection-based reduced order modeling and we assume, as it is common practice in this case, that the eigenvalue problem is affine parameter dependent. Also, we follow the offline-online paradigm for our approach. In the offline stage, the discretized parametric eigenvalue problem is solved for some sample parameters using any numerical methods like FDM, FEM, FVM, and a reduced space is constructed using those solutions. Then projecting the problem into this reduced space leads to a reduced problem. In the online stage, for a given value of the parameter, we solve the reduced problem and use the reduced solution as a proxy to the corresponding solution in the original discrete space. 

In our study we use the finite element method for the solution ot the discretized eigenvalue problem in the offline stage. The variational form of the parametric PDE we are considering, reads as follows: given $\pmb{\mu}\in\mathcal{P}$ find $\lambda(\pmb{\mu}) \in \mathbb{R}$ and $u(\pmb{\mu})\neq 0\in V$ such that
\begin{align}\label{var}
    a(u,v;\pmb{\mu})=\lambda(\pmb{\mu}) b(u,v;\pmb{\mu})\quad \forall v\in V,
\end{align}
where $a,b:V\times V: \to \mathbb{R}$ are symmetric bilinear forms defined on a Hilbert space $V$. The parameter set $\mathcal{P} \subset \mathbb{R}^p$ is closed and bounded.
Let $V_h$ be a finite dimensional subspace of $V$. Then the discretized problem corresponding to the variational form \eqref{var} reads: given $\pmb{\mu}\in\mathcal{P}$ find $\lambda_h(\pmb{\mu}) \in \mathbb{R}$ and $u_h(\pmb{\mu})\neq 0\in V$ such that
\begin{align}\label{varh}
    a(u_h,v_h;\pmb{\mu})=\lambda_h(\pmb{\mu}) b(u_h,v_h;\pmb{\mu})\quad \forall v_h\in V_h.
\end{align} 
This reduces to a generalized eigenvalue problem 
\begin{align}\label{mat_hf}
    A_h^{\pmb{\mu}} \pmb{u}_h=\lambda_h B_h^{\pmb{\mu}} \pmb{u}_h
\end{align}
where $A_h^{\pmb{\mu}}(i,j)=a(\varphi_j,\varphi_i;\pmb{\mu})$ and $B_h^{\pmb{\mu}}(i,j)=b(\varphi_j,\varphi_i;\pmb{\mu})$ for $i,j=1,2,\dots, N_h$ and $\{\varphi_1,\dots,\varphi_{N_h} \}$ is a basis of $V_h$. The vector $\pmb{u}_h$ contains the nodal values of the solution $u_h$.
In addition to the symmetry of the bilinear forms $a,b$, we assume that they are continuous, that $a$ is $V$-elliptic and that $b$ is equivalent to a scalar product in $L^2(\Omega)$. For a more detailed discussion on the assumptions on the bilinear forms we refer the reader to~\cite{Boffi10, Moataz}. The bilinear forms for an affinely parameter dependent eigenvalue problem are as follows:
\begin{align*}
    a(w,v;\pmb{\mu})=\sum\limits_{l=1}^{n_a} \theta_l^a (\pmb{\mu} )a_l(w,v),\\ \quad b(w,v;\pmb{\mu})=\sum\limits_{m=1}^{n_b} \theta_m^b (\pmb{\mu} )b_m(w,v),
\end{align*}
where $a_l,b_m:V \times V \to \mathbb{R}$ are bilinear forms that are free from the parameters. In this case the matrix form  \eqref{mat_hf} of the eigenvalue problem is given by:
\begin{align}\label{mat_hf_afn}
     \sum\limits_{l=1}^{n_a} \theta_l^a (\pmb{\mu} ) A_l^h\pmb{u}_h=\lambda_h(\pmb{\mu}) \sum\limits_{m=1}^{n_b} \theta_m^b (\pmb{\mu} ) B_m^h \pmb{u}_h
\end{align}
where $A_l^h(i,j)=a_l(\varphi_j,\varphi_i)$ and $B_m^h(i,j)=b_m(\varphi_j,\varphi_i)$.
\section{Reduced order modeling}\label{sec:rom}
Let us suppose that the eigenvalue problem has been solved for parameters $\pmb{\mu}_1,\pmb{\mu}_2,\dots, \pmb{\mu}_{n_s}$. Using these solutions we find an orthonormal basis for the reduced space using the POD approach. In the POD approach, first we form the snapshot matrix $S$ using eigenvectors calculated at the sample parameters. The choice of the eigenvectors is very crucial for getting good approximation of the eigenvalues and eigenvectors from the reduced system. We have tested possible choices of different combinations of eigenvectors for the snapshot matrix in~\cite{WithGopal}. The snapshot matrix has size $N_h\times n_k$ with $n_k=n_e.n_s$, where $n_e$ denotes the number of eigenvectors we are considering in the snapshot matrix at one particular sample point.  Using singular value decomposition on the snapshot matrix $S$, we get two orthogonal matrices $Q\in \mathbb{R}^{N_h\times N_h}$ and $R\in \mathbb{R}^{n_k\times n_k}$ such that
$$S=Q \Sigma R^\top,$$
where $\Sigma$ is a diagonal matrix of size $N_h\times n_k$ containing the singular values. If $r$ is the rank of $S$, then the first $r$ diagonal elements of $\Sigma$ are non-zero. The column vectors of the matrix $Q$ are called the left singular vectors of the matrix $S$ and the columns of $R$ are called the right singular vectors of the matrix $S$. The first $N$ dominating left singular vectors are treated as the basis of the reduced space. 
For a given tolerance $\epsilon_{tol}$, we will consider $N$ as the smallest integer satisfying the criterion
\begin{equation}\label{criterion}
\frac{\sum\limits_{i = 1}^N \sigma_{i}^2}{\sum\limits_{i = 1}^r \sigma_{i}^2} \geq 1 - \epsilon_{tol}
\end{equation}
In our numerical computation we use $\epsilon_{tol}=10^{-8}$. The criterion \eqref{criterion}, follows from the Schmidt--Eckart--Young theorem which states that if a matrix $M$ of rank $r$ is approximated by $M_k=\sum\limits_{i=1}^k \sigma_i \pmb{\xi}_i \pmb{\eta}_i^T$ using the first $k$ left singular vectors $\{\pmb{\xi}_i\}_{i=1}^k$ and right eigenvectors $\{\pmb{\eta}_i\}_{i=1}^k$ of $M$, then the square of the approximation error ($\|M-M_k\|_F^2$) in the Frobenius norm is equal to the sum of the squares of last $r-k$ singular values.

The left singular vectors can be obtained using the eigenvectors of the Gram matrix $S^\top S$. The choice of the reduced basis consisting of left singular vectors is optimal in the sense that it minimizes the sum of the squares of the error of each snapshot $\pmb{u}_i$ ($i^{th}$ column of $S$) and its projection onto any $N$-dimensional subspace~\cite{Quarteronietal16}.

Let $\pmb{\zeta}_1,\dots, \pmb{\zeta}_N$ be the first $N$ left singular vectors (i.e. first $N$ columns of Q) of the snapshot matrix then the reduced space is defined as $V_N:=span\{\zeta_1,\zeta_2,\dots,\zeta_N\}$ where $$\zeta_j=\sum\limits_{i=1}^{N_h} \zeta_j^{(i)}\varphi_i, \quad j=1,\dots, N$$ 
and $$\pmb{\zeta}_k=(\zeta^{(1)}_k,\dots,\zeta^{(N_h)}_k)^T \quad k=1,\dots, N.$$ Then, projecting the variational form into the reduced space, we obtain the reduced eigenvalue problem in the variational form: given $\pmb{\mu}\in \mathcal{P}$, find $\lambda_N\in \mathbb{R}$, $u_N(\pmb{\mu})\neq 0 \in V_N$ such that 
\begin{align}\label{red_var}
    a(u_N,v_N;\pmb{\mu})=\lambda_N(\pmb{\mu}) b(u_N,v_N;\pmb{\mu})\quad \forall v_N\in V_N.
\end{align}
Note that from $V_N \subseteq V_h\subseteq V$ we obain $\lambda \leq \lambda_h \leq \lambda_N$ as the eigenvalues are characterized by the minmax property of the Rayleigh quotient.
Let $u_N=\sum\limits_{j=1}^N u^{(j)}_N\zeta_j$, substituting this expression of $u_N$ in \eqref{red_var} and choosing $v_N=\zeta_i,~i=1,\dots,N$, we have
\begin{align*}
    A_N^{\pmb{\mu}} \pmb{u}_N=\lambda_N(\pmb{\mu})B_N^{\pmb{\mu}} \pmb{u}_N,
\end{align*}
where $$A_N^{\pmb{\mu}}(i,j)=a(\zeta_j,~\zeta_i;\pmb{\mu}),~~B_N^{\pmb{\mu}}(i,j)=b(\zeta_j,\zeta_i;\pmb{\mu})$$ and  $\pmb{u}_N=(u^{(1)}_N,\dots,u^{(N)}_N)^\top$ is the vector containing the coefficients of the reduced solution $u_N$. 
Assuming the affine parameter dependent condition, the matrix form of the reduced system will be
\begin{align}\label{mat_red_afn}
     \sum\limits_{l=1}^{n_a} \theta_l^a (\pmb{\mu} ) A_l^N\pmb{u}_N=\lambda_N(\pmb{\mu})\sum\limits_{m=1}^{n_b} \theta_m^b (\pmb{\mu} ) B_m^N \pmb{u}_N,
\end{align}
where $A_l^N(i,j)=a_l(\zeta_j,\zeta_i)$ and $B_m^N(i,j)=b_m(\zeta_j,\zeta_i)$. Using the expression of $\zeta_i$ and $\zeta_j$, we have 
$$a_l(\zeta_j,\zeta_i)=\sum\limits_{r=1}^{N_h}\sum\limits_{s=1}^{N_h}\zeta^{(s)}_ja_l(\varphi_j,\varphi_i)\zeta^{(r)}_i\quad 1\leq i,j\leq N$$
and $$b_m(\zeta_j,\zeta_i)=\sum\limits_{r=1}^{N_h}\sum\limits_{s=1}^{N_h}\zeta^{(s)}_jb_m(\varphi_j,\varphi_i)\zeta^{(r)}_i\quad 1\leq i,j\leq N.$$
Let $\mathbb{V}=[\pmb{\zeta}_1 \dots \pmb{\zeta}_N]\in \mathbb{R}^{N_h\times N}$ be the transformation matrix with
$$\mathbb{V}(i,j)=\zeta^{(i)}_j, \quad 1\leq i \leq N_h,~ 1\leq j\leq N .$$
Thus we can write $$A_l^N=\mathbb{V}^\top A_l^h\mathbb{V},\quad B_m^N=\mathbb{V}^\top B_m^h\mathbb{V}.$$
Thus the matrix from \eqref{mat_red_afn} of the reduced problem for the affine dependent problem reduces to:
\begin{equation}\label{red_sys}
      \sum\limits_{l=1}^{n_a} \theta_l^a (\pmb{\mu} ) \mathbb{V}^\top A_l^h \mathbb{V}\pmb{u}_N=\lambda_N(\pmb{\mu})\sum\limits_{m=1}^{n_b} \theta_m^b (\pmb{\mu} ) \mathbb{V}^\top B_m^h\mathbb{V} \pmb{u}_N.
\end{equation}
In the offline stage when we solve the problem for the sample parameters $\pmb{\mu}_1,\dots,\pmb{\mu}_{n_s}$, we have to calculate the parameter independent matrices $A^h_l,~B^h_m$ and the transformation matrix $\mathbb{V}$. Thus the matrices $\mathbb{V}^\top A_l^h\mathbb{V}$ and $\mathbb{V}^\top B_m^h\mathbb{V}$ are already known from the offline phase. In the online phase, when we are given a new parameter $\pmb{\mu}_{\ast}$, we can obtain the reduced system \eqref{red_sys} just by evaluating the parameter dependent functions $\theta_l^a$ and $\theta_m^b$. It is remarkable that solving the reduced system is in general much faster than the original problem if its dimension $N$ is much smaller than $N_h$, the size of the FEM problem. In practice we have that $N\ll N_h$.

Once we have the reduced solution vector $\pmb{u}_N$, we get the reduced solution by using its expression $$u_N=\sum\limits_{j=1}^N u^{(j)}_N\zeta_j=\sum\limits_{j=1}^N \sum\limits_{i}^{N_h} u^{(j)}_N \zeta_j^{(i)}\varphi_i. $$
Thus the reduced solution vector $\pmb{u}_N$ and the corresponding high fidelity (FEM) solution vector $\pmb{u}_h$ are related by 
$$\pmb{u}_h=\mathbb{V}\pmb{u}_N, \quad \pmb{u}_N=\mathbb{V}^\top\pmb{u}_h. $$

\section{Numerical results}\label{sec:nr}
In this section we present the discretization of two parametric eigenvalue problems using four different types of sampling techniques, 
namely random sampling, Latin hypercube sampling, Smolyak sampling, and uniform sampling. Here we compared the first three eigenvalues obtained using FEM and ROM with reduced basis obtained from the snapshots corresponding to the sampling points. In order to quantify the error between the FEM eigenvalues and reduced eigenvalues, we have reported the relative error in the corresponding tables.
\subsection{An eigenvalue problem with two parameters}
\begin{figure}
     \begin{subfigure}{0.22\textwidth}
        \includegraphics[height=3.5cm,width=4.2cm]{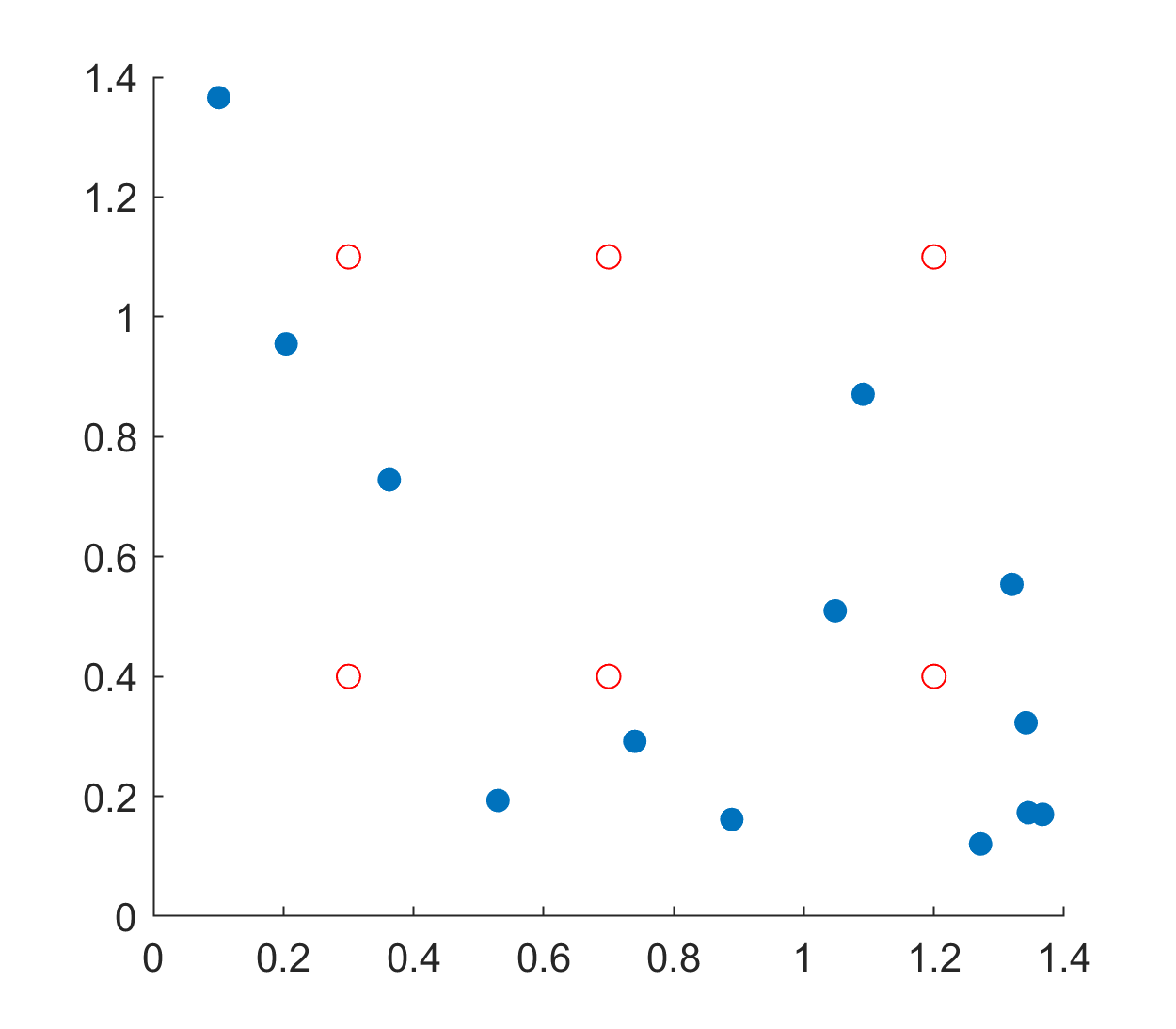}
       \caption{Random Sampling}
     \end{subfigure}
     \begin{subfigure}{0.22\textwidth}
          \includegraphics[height=3.5cm,width=4.2cm]{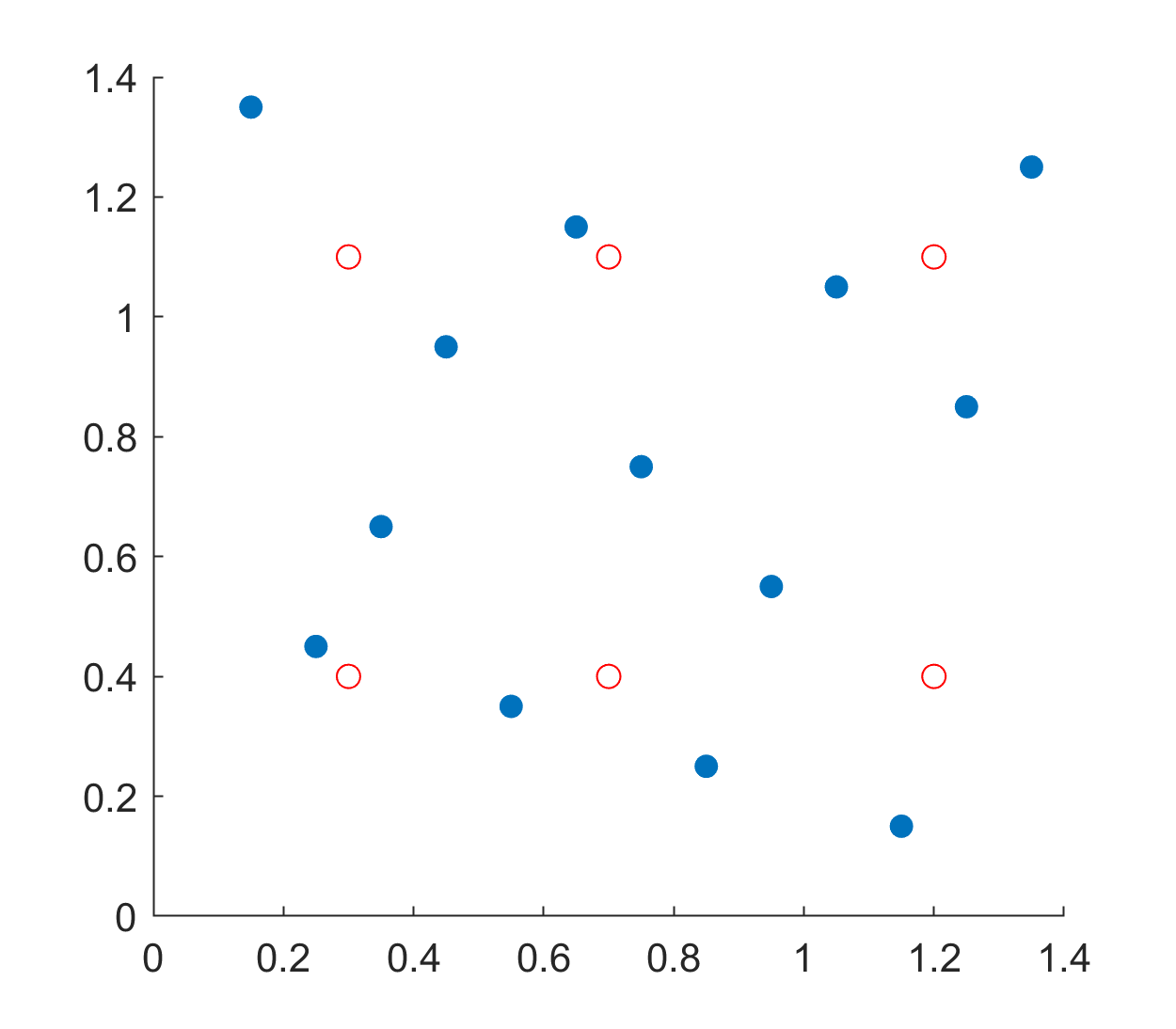}
         \caption{LHS Sampling}
     \end{subfigure}
     \begin{subfigure}{0.22\textwidth}
          \includegraphics[height=3.5cm,width=4.2cm]{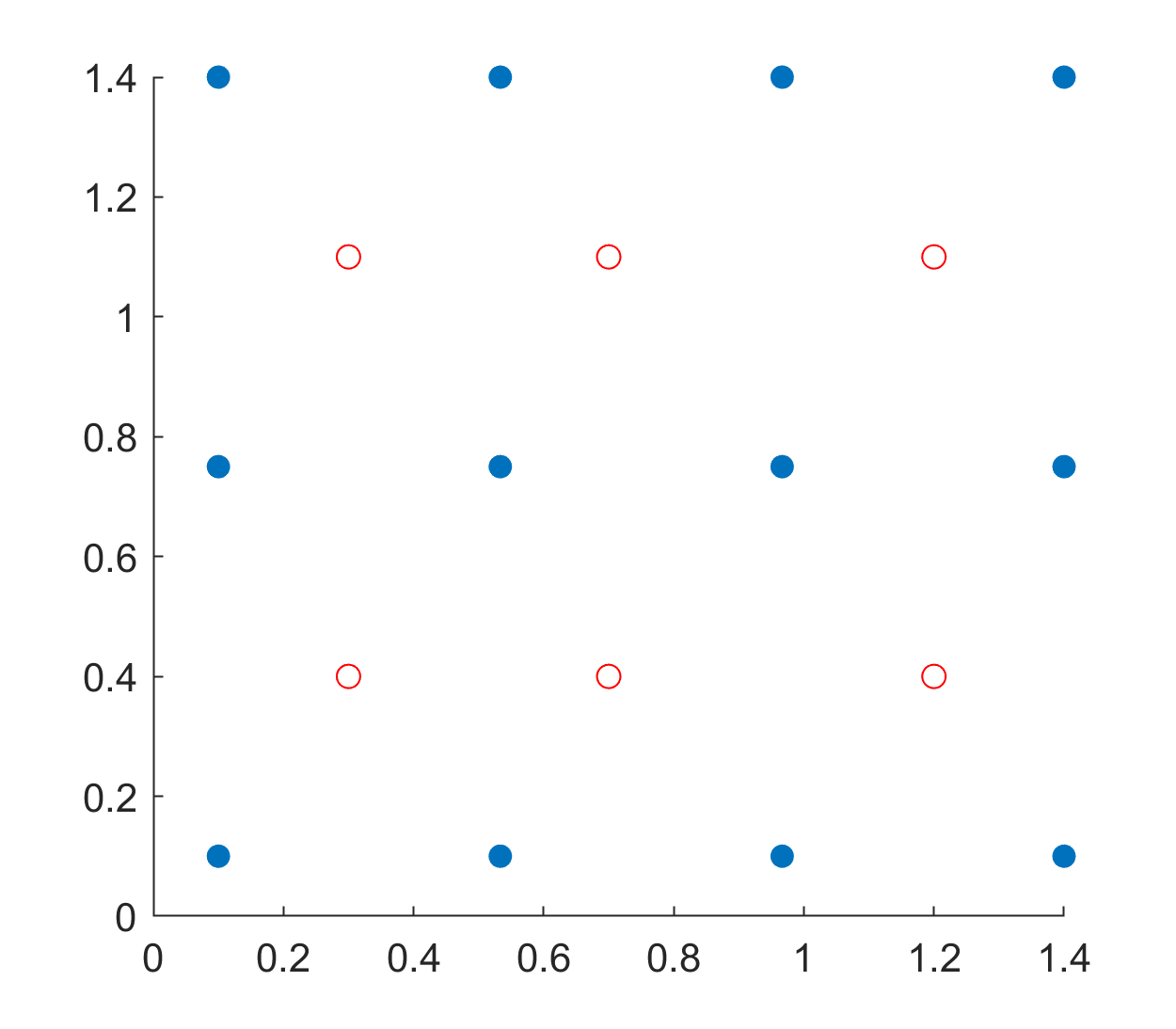}
         \caption{Uniform Sampling}
     \end{subfigure}
      \begin{subfigure}{0.22\textwidth}
          \includegraphics[height=3.5cm,width=4.2cm]{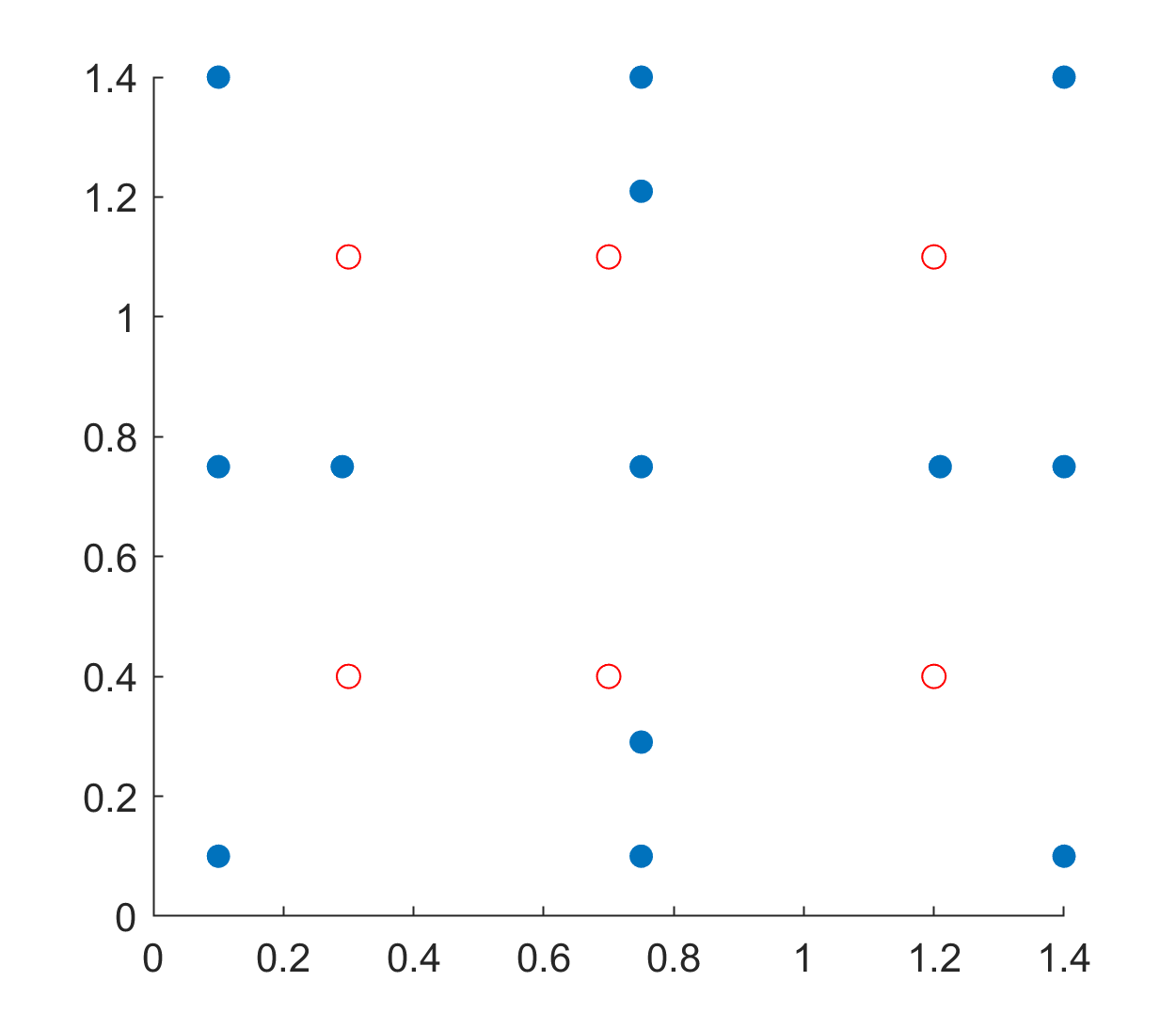}
         \caption{Smolyak Sampling}
     \end{subfigure}
  \caption{Different sampling of $[0.1,1.4]^2$ with 13 points (in blue) and 6 test points (in red).}
        \label{sample_2d}
  \end{figure}

Let us consider the following eigenvalue problem
 \begin{equation}
\label{eq:PDE_example}
\left\{
\begin{array}{ll}
     -\operatorname{div}(A(\pmb{\pmb{\mu}})\nabla u(\pmb{\pmb{\mu}}))=\lambda(\pmb{\pmb{\mu}}) u(\pmb{\pmb{\mu}})&  \textrm{ in }\Omega=(0,1)^2\\
     u(\pmb{\pmb{\mu}})=0& \textrm{ on }\partial\Omega,
\end{array}
\right.
\end{equation}
where the diffusion $A(\pmb{\mu})\in\mathbb R^{2\times 2}$ is given by the matrix
\begin{equation}\label{Amu}
    A(\pmb{\mu})= 
    \begin{pmatrix}
        \frac{1}{\mu_1^2} &\frac{0.7}{\mu_2}\\
        \frac{0.7}{\mu_2} &\frac{1}{\mu_2^2}
    \end{pmatrix},
\end{equation}
with $\pmb{\mu}=(\mu_1,\mu_2)\in \mathcal{P}\subset \mathbb{R}^2.$
The problem is symmetric and the parameter space $\mathcal{P}$ is chosen in such a way that the matrix is positive definite. The matrix is positive definite for any nonzero value of $\mu_2$ and $\mu_1\in (-1.42,1.42)\setminus\{0\}$. For the numerical purpose, we choose the parameter space to be $\mathcal{P}=[0.1,1.4]^2$. In Figure~\ref{sample_2d}, we displayed the 13 sample points using different techniques from the parameter domain $\mathcal{M}$, and also we displayed six test points in red that are used to compare the results of FEM and ROM. In order to discretize the domain we have used the MATLAB command \texttt{initmesh} with maximum mesh size $h=0.01.$

In Table~\ref{table:u1_p2} we report the first eigenvalues at the six test points obtained by FEM and ROM corresponding to different sampling; we include also the relative error between them. We formed the snapshot matrix by considering only the first eigenvector at the sample points. The number of POD  basis obtained using the criterion is mentioned in the bracket of the first column.
From the table, one can see that the eigenvalues corresponding to LHS sampling are better than the uniform sampling.
We report the second eigenvalues in Table~\ref{table:u12_p2} and the third eigenvalues in Table~\ref{table:u13_p2}. As anticipated in~\cite{WithGopal}, we have used the first and the second eigenvectors in the snapshot matrix for calculating second eigenvalues and  first, second, and third eigenvalues for calculating third eigenvalues, respectively. Also in these cases, we can see that the results corresponding to the LHS sampling are better than the uniform sampling, while other choices are giving comparable results.


\begin{table}
 	 	\centering
 	\begin{tabular}{|c|c|c|c|c|c|c|c|} 
 		\hline
 	& 
		  {\begin{tabular}[c]{@{}c@{}} $\mu$ \end{tabular}} &
 		 {\begin{tabular}[c]{@{}c@{}}  FEM \end{tabular}} & 
 		{\begin{tabular}[c]{@{}c@{}}  ROM  \end{tabular}}&
 		{\begin{tabular}[c]{@{}c@{}}   Error  \end{tabular}}	\\
  \hline
  \parbox[t]{2mm}{\multirow{6}{*}{\rotatebox[origin=c]{90}{Random  (11)}}} 
 &(0.3,0.4)&169.9519383    &169.9552761        &$1.9\times 10^{-5}$\\
  &(0.3,1.1)&117.5380721    &117.5381817        &$9.3\times 10^{-7}$\\
      &(0.7,0.4)&78.6821102     &78.6827105     &$7.6\times 10^{-6}$\\
      &(0.7,1.1)&27.0986848     &27.1034374     &$1.7\times 10^{-4}$ \\
      &(1.2,0.3)&64.2706410     &64.2713395     &$1.0\times 10^{-5}$\\
      &(1.2,1.1)&12.1778849     &12.1801075     &$1.8\times 10^{-4}$\\
  \hline
  \parbox[t]{2mm}{\multirow{6}{*}{\rotatebox[origin=c]{90}{LHS (11)}}} 
  &(0.3,0.4)&169.9519383     &169.9521614       &$1.3\times v{-6}$\\
  &(0.3,1.1)&117.5380721    &117.5380728     &$6.1\times  10^{-9}$\\
      &(0.7,0.4)&78.6821102     &78.6821284     &$2.3\times  10^{-7}$\\
      &(0.7,1.1)&27.0986848     &27.0988903      &$7.5\times  10^{-6}$ \\
      &(1.2,0.3)&64.2706410     &64.2744779      &$5.9\times  10^{-5}$\\
      &(1.2,1.1)&12.1778849     &12.1779026     &$1.4\times  10^{-6}$\\
  \hline
  \parbox[t]{2mm}{\multirow{6}{*}{\rotatebox[origin=c]{90}{Smolyak (11)}}} 
 &(0.3,0.4)&169.9519383   &169.9525321           &$3.4\times  10^{-6}$\\
  &(0.3,1.1)&117.5380721    &117.5382384         &$1.4\times  10^{-6}$\\
      &(0.7,0.4)&78.6821102     &78.6825511      &$5.6\times  10^{-6}$\\
      &(0.7,1.1)&27.0986848     &27.0986968      &$4.4\times  10^{-7}$ \\
      &(1.2,0.3)&64.2706410     &64.2754208      &$7.4\times  10^{-5}$\\
      &(1.2,1.1)&12.1778849     &12.1838340      &$4.8\times  10^{-4}$\\
      \hline
      \parbox[t]{2mm}{\multirow{6}{*}{\rotatebox[origin=c]{90}{Uniform  (11)}}} 
  &(0.3,0.4)&169.9519383    &169.9520650      &$7.4\times 10^{-7}$\\
  &(0.3,1.1)&117.5380721    &117.5384261     &$3.0\times 10^{-6}$\\
      &(0.7,0.4)&78.6821102     &78.6842588       &$2.7\times 10^{-5}$\\
      &(0.7,1.1)&27.0986848     &27.0987477       &$2.3\times 10^{-6}$ \\
          &(1.2,0.3)&64.2706410     &64.2837808       &$2.0\times 10^{-4}$\\
          &(1.2,1.1)&12.1778849     &12.1823100      &$3.6\times 10^{-4}$\\
      \hline
\end{tabular}
\caption{Comparison of first eigenvalues of \eqref{eq:PDE_example}, obtained using FEM and ROM: the snapshot matrix contains $u_1$ at sample points.} 
\label{table:u1_p2}
 \end{table}

\begin{table}
 	 	\centering
 	\begin{tabular}{|c|c|c|c|c|c|c|c|} 
 		\hline
 	& 
		  {\begin{tabular}[c]{@{}c@{}} $\mu$ \end{tabular}} &
 		 {\begin{tabular}[c]{@{}c@{}}  FEM \end{tabular}} & 
 		{\begin{tabular}[c]{@{}c@{}}  ROM  \end{tabular}}&
 		{\begin{tabular}[c]{@{}c@{}}   Error  \end{tabular}}	\\
  \hline
  \parbox[t]{2mm}{\multirow{3}{*}{\rotatebox[origin=c]{90}{Random(24)}}} 
&(0.3,0.4)&348.5059152     &348.5261134     &$5.7\times 10^{-5}$\\
&(0.3,1.1)&141.1559804     &141.1564666     &$3.4\times 10^{-6}$\\
&(0.7,0.4)&128.5084550     &128.5160517     &$5.9\times 10^{-5}$\\
&(0.7,1.1)&47.2728621      &47.2976190      &$5.2\times 10^{-4}$ \\
&(1.2,0.3)&71.9749059      &71.9783770      &$4.8\times 10^{-5}$\\
&(1.2,1.1)&20.8445384      &20.8699567      &$1.2\times 10^{-3}$\\
  \hline
\parbox[t]{2mm}{\multirow{3}{*}{\rotatebox[origin=c]{90}{LHS(22)}}} 
&(0.3,0.4)&348.5059152     &348.5072206       &$3.7\times 10^{-6}$\\
&(0.3,1.1)&141.1559804     &141.1559835      &$2.1\times 10^{-8}$\\
&(0.7,0.4)&128.5084550     &128.5087958     &$2.6\times 10^{-6}$\\
&(0.7,1.1)&47.2728621      &47.2736609        &$1.6\times 10^{-5}$ \\
&(1.2,0.3)&71.9749059      &71.9885937       &$1.9\times 10^{-4}$\\
&(1.2,1.1)&20.8445384      &20.8446366      &$4.7\times 10^{-6}$\\
\hline
\parbox[t]{2mm}{\multirow{3}{*}{\rotatebox[origin=c]{90}{Smolyak(23)}}} 
&(0.3,0.4)&348.5059152     &348.5241038      &$5.2\times 10^{-5}$\\
&(0.3,1.1)&141.1559804     &141.1586502      &$1.8\times 10^{-5}$\\
&(0.7,0.4)&128.5084550     &128.5109413      &$1.9\times 10^{-5}$\\
&(0.7,1.1)&47.2728621      &47.2729549       &$1.9\times 10^{-6}$ \\
&(1.2,0.3)&71.9749059      &72.0088014       &$4.7\times 10^{-4}$\\
&(1.2,1.1)&20.8445384      &20.8805726       &$1.7\times 10^{-3}$\\
\hline  
\parbox[t]{2mm}{\multirow{3}{*}{\rotatebox[origin=c]{90}{Smolyak(22)}}} 
&(0.3,0.4)&348.5059152     &348.5101071       &$1.2\times 10^{-4}$\\
&(0.3,1.1)&141.1559804     &141.1570619        &$7.6\times 10^{-6}$\\
&(0.7,0.4)&128.5084550     &128.5712840      &$4.8\times 10^{-4}$\\
&(0.7,1.1)&47.2728621      &47.2730675        &$4.3\times 10^{-6}$ \\
&(1.2,0.3)&71.9749059      &72.1402955        &$2.3\times 10^{-3}$\\
&(1.2,1.1)&20.8445384      &20.9061857       &$3.0\times 10^{-3}$\\
      \hline  
\end{tabular}
\caption{Comparison of second eigenvalues of \eqref{eq:PDE_example}, obtained using FEM and ROM: the snapshot matrix contains $u_1,u_2$ at sample points.} 
\label{table:u12_p2}
 \end{table}

\begin{table}
 	 	\centering
 	\begin{tabular}{|c|c|c|c|c|c|c|c|} 
 		\hline
 	& 
		  {\begin{tabular}[c]{@{}c@{}} $\mu$ \end{tabular}} &
 		 {\begin{tabular}[c]{@{}c@{}}  FEM \end{tabular}} & 
 		{\begin{tabular}[c]{@{}c@{}}  ROM  \end{tabular}}&
 		{\begin{tabular}[c]{@{}c@{}}   Error  \end{tabular}}	\\
  \hline
  \parbox[t]{2mm}{\multirow{3}{*}{\rotatebox[origin=c]{90}{Random (37)}}} 
  &(0.3,0.4)&501.1398704    &501.3300166    &$3.7\times 10^{-4}$\\
  &(0.3,1.1)&180.5284070    &180.5293306    &$5.1\times 10^{-6}$\\
      &(0.7,0.4)&206.5455501    &206.5521503    &$3.1\times 10^{-5}$\\
      &(0.7,1.1)&77.7273406     &77.7538331     &$3.4 \times 10^{-4}$ \\
      &(1.2,0.3)&84.6527187     &84.6566887     &$4.6\times 10^{-5}$\\
      &(1.2,1.1)&29.8899408     &29.9277991     &$1.3\times 10^{-3}$\\
  \hline
  \parbox[t]{2mm}{\multirow{3}{*}{\rotatebox[origin=c]{90}{LHS (37)}}} 
 &(0.3,0.4)&501.1398704    &501.1592711      &$3.8\times 10^{-5}$\\
  &(0.3,1.1)&180.5284070   &180.5284112      &$2.3\times 10^{-8}$\\
      &(0.7,0.4)&206.5455501    &206.5458945     &$1.6\times 10^{-6}$\\
      &(0.7,1.1)&77.7273406     &77.7304915      &$4.0 \times 10^{-5}$ \\
      &(1.2,0.3)&84.6527187     &84.6737108      &$2.4\times 10^{-4}$\\
      &(1.2,1.1)&29.8899408     &29.8901000      &$5.3\times 10^{-6}$\\
  \hline
  \parbox[t]{2mm}{\multirow{3}{*}{\rotatebox[origin=c]{90}{Smolyak (38)}}} 
  &(0.3,0.4)&501.1398704   &501.1503121     &$2.0\times 10^{-5}$\\
  &(0.3,1.1)&180.5284070   &180.5295786      &$6.4\times 10^{-6}$\\
      &(0.7,0.4)&206.5455501    &206.5552724     &$4.7\times 10^{-5}$\\
      &(0.7,1.1)&77.7273406     &77.7275723      &$2.9 \times 10^{-6}$ \\
      &(1.2,0.3)&84.6527187     &84.7392639      &$1.0\times 10^{-3}$\\
      &(1.2,1.1)&29.8899408     &29.9374827      &$1.6\times 10^{-3}$\\
  \hline
  \parbox[t]{2mm}{\multirow{3}{*}{\rotatebox[origin=c]{90}{Uniform (36)}}} 
  &(0.3,0.4)&501.1398704    &501.1458402       &$1.1\times 10^{-5}$\\
  &(0.3,1.1)&180.5284070   &180.5303306     &$1.0\times 10^{-5}$\\
      &(0.7,0.4)&206.5455501    &206.7809186   &$1.1\times 10^{-3}$\\
      &(0.7,1.1)&77.7273406     &77.7383395     &$1.4 \times 10^{-4}$ \\
      &(1.2,0.3)&84.6527187     &85.0786560    &$5.0\times 10^{-3}$\\
      &(1.2,1.1)&29.8899408     &29.9486247    &$2.0\times 10^{-3}$\\
  \hline
\end{tabular}
\caption{Comparison of third eigenvalues of \eqref{eq:PDE_example}, obtained using FEM and ROM: the snapshot matrix contains $u_1,u_2,u_3$ at sample points.} 
\label{table:u13_p2}
 \end{table}

\subsection{An eigenvalue problem with three parameters}
\begin{figure}
     \begin{subfigure}{0.22\textwidth}
        \includegraphics[height=3cm,width=4.2cm]{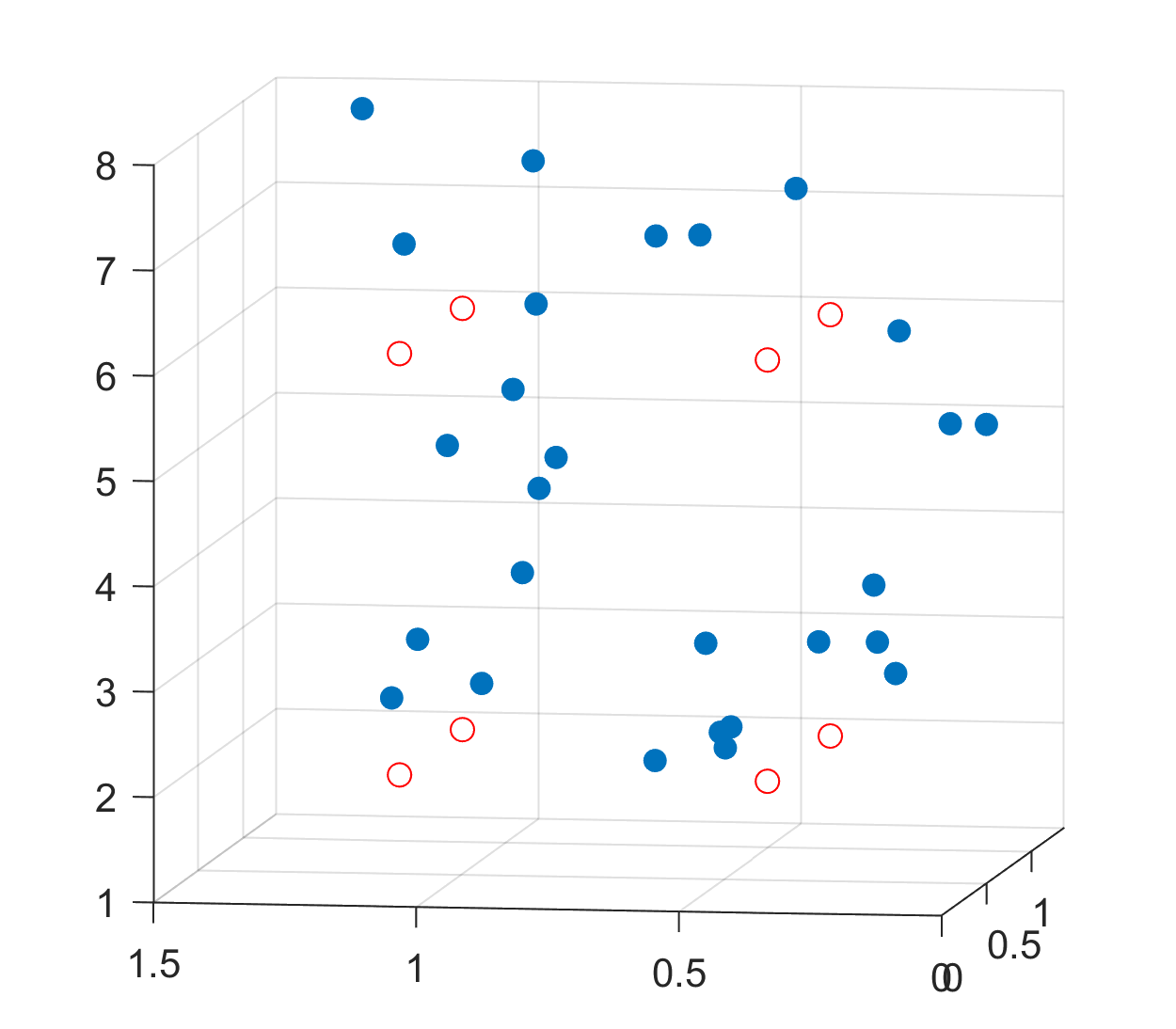}
       \caption{Random Sampling}
     \end{subfigure}
     \begin{subfigure}{0.22\textwidth}
          \includegraphics[height=3cm,width=4.2cm]{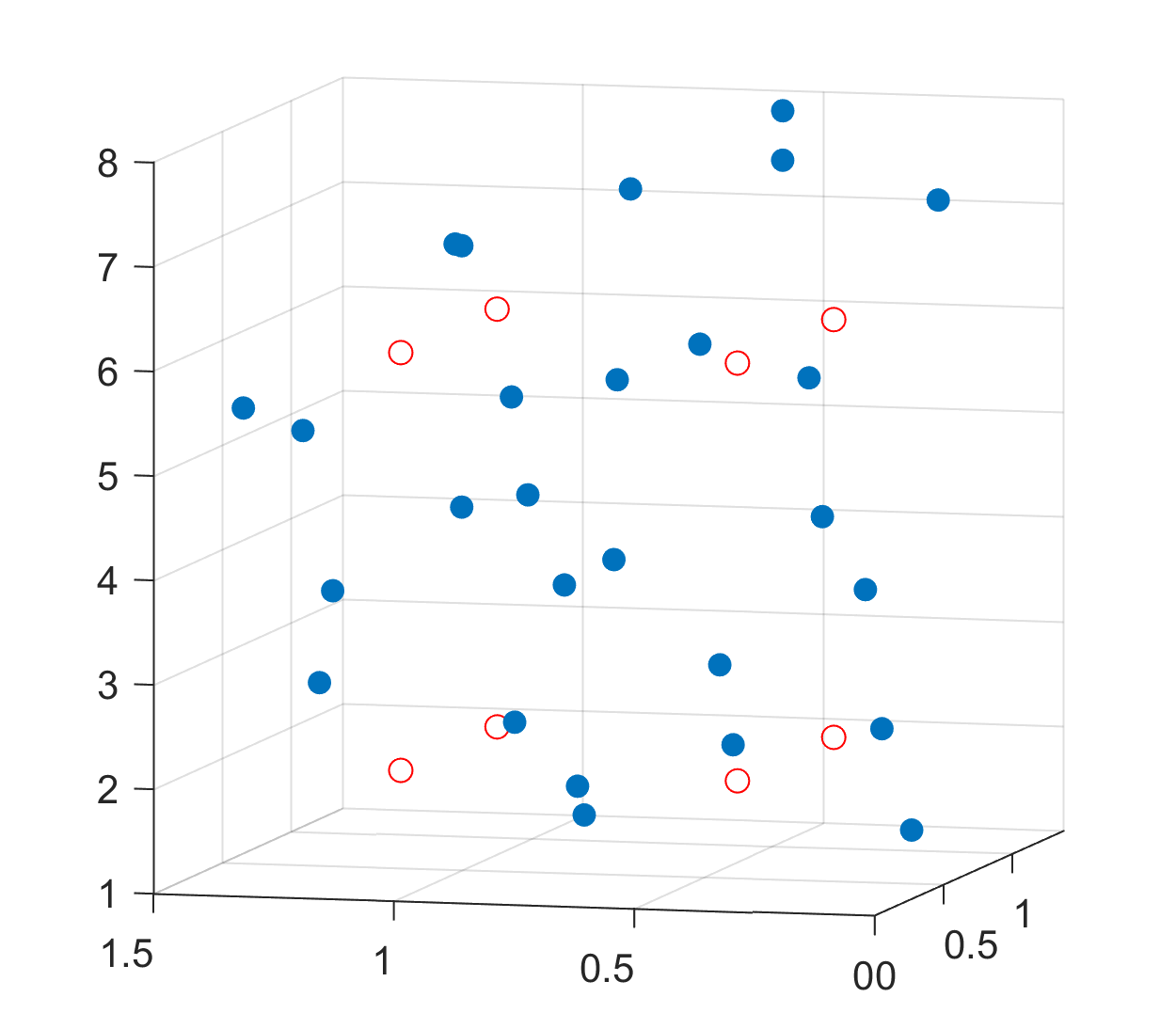}
         \caption{LHS Sampling}
     \end{subfigure}
     \begin{subfigure}{0.22\textwidth}
          \includegraphics[height=3cm,width=4.2cm]{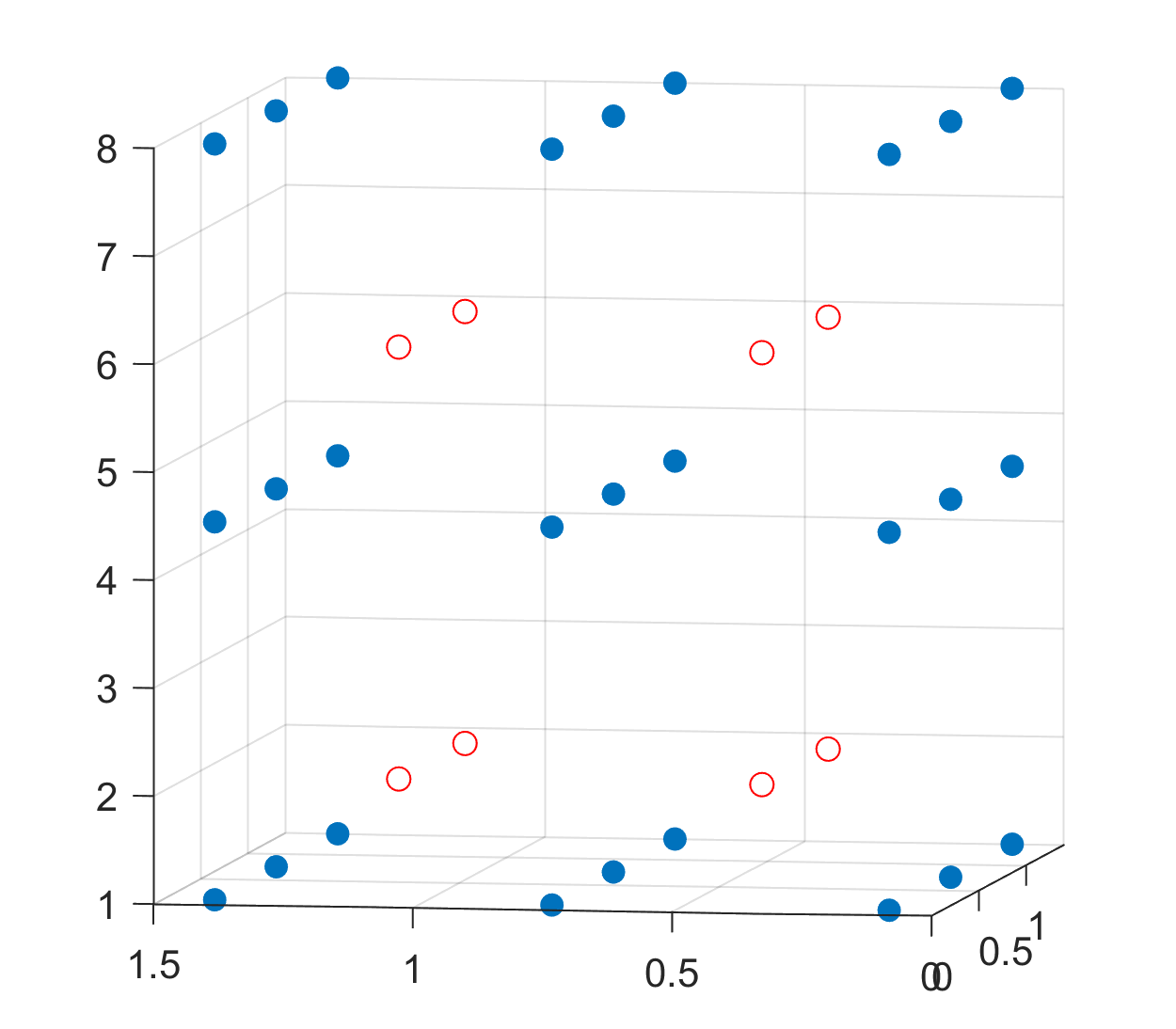}
         \caption{Uniform Sampling}
     \end{subfigure}
      \begin{subfigure}{0.22\textwidth}
          \includegraphics[height=3cm,width=4.2cm]{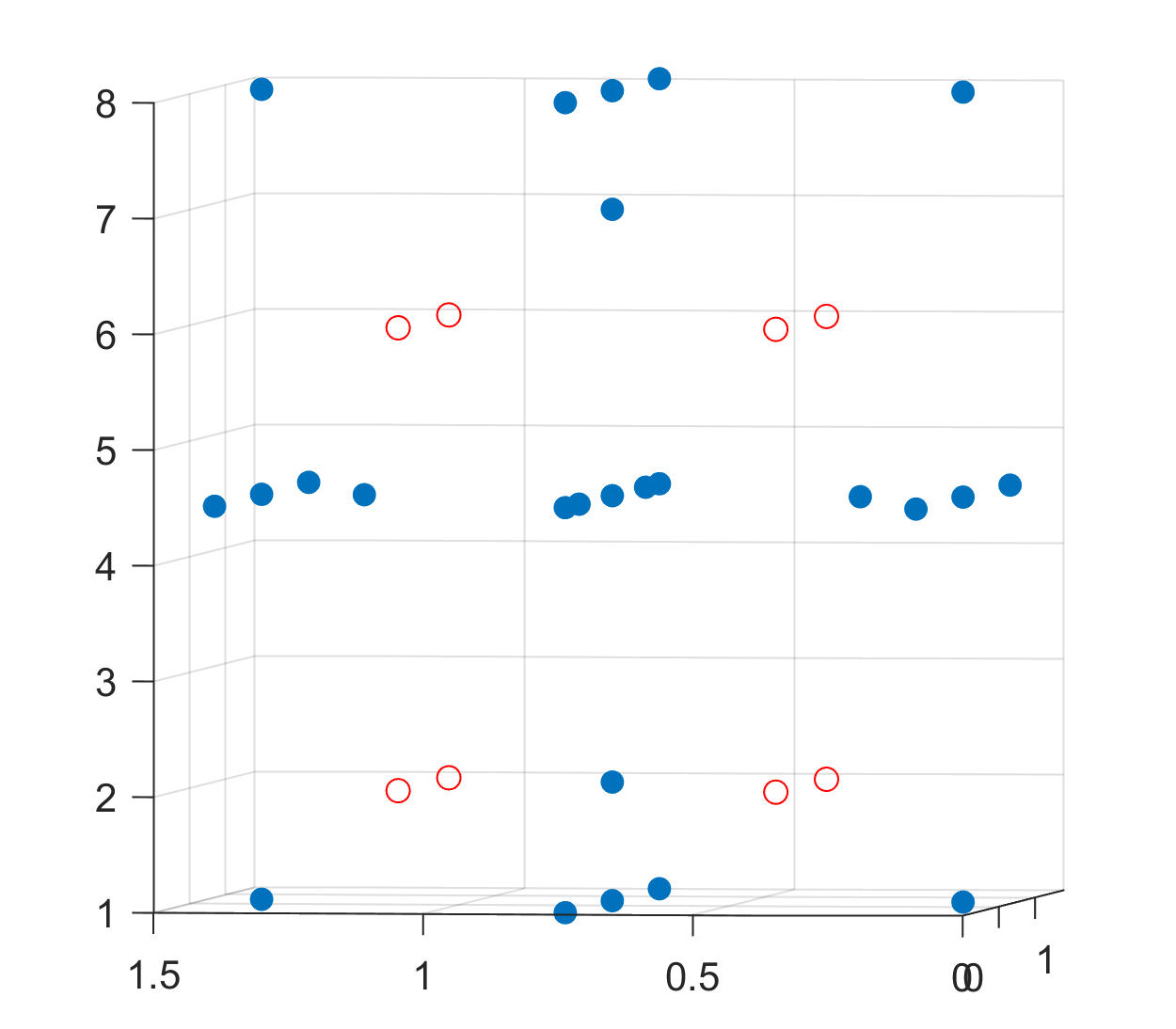}
         \caption{Smolyak Sampling}
     \end{subfigure}
  \caption{Different sampling of $[0.4,1.4]^2\times [1,8]$ with 27 points (in blue). Also plotted 8 test points (in red).}
        \label{sample_3d}
  \end{figure}

Let us consider the following eigenvalue problem with three parameters:
 \begin{equation}
\label{eq:PDE_example1}
\left\{
\aligned
&-\operatorname{div}(A(\pmb{\mu})\nabla u(\pmb{\mu}))+\frac{\mu_3^2}{2}(x^2+y^2)u(\pmb{\mu})=\lambda(\pmb{\mu})u(\pmb{\mu})&& \textrm{in } \Omega\\
&u(\pmb{\mu})=0&& \textrm{on }\partial\Omega,
\endaligned
\right.
\end{equation}
where the diffusion $A(\pmb{\mu})\in\mathbb R^{2\times 2}$ is defined in \eqref{Amu}
with $\pmb{\mu}=(\mu_1,\mu_2,\mu_3)\in \mathcal{P}\subset \mathbb{R}^3.$
 For the numerical purpose, we choose the parameter space to be $\mathcal{P}=[0.1,1.4]\times [0.1,1.4] \times [1,8]$. In Figure~\ref{sample_3d} we display the 27 sample points using different techniques from the parameter domain $\mathcal{M}$, and we display also eight test points in red that are used to compare the results of FEM and ROM.

In Table~\ref{table:u1_p3} we report the first eigenvalues at the test points obtained by FEM and the reduced order model corresponding to different samplings. In the snapshot matrix, we have included only the first eigenvector at the sample points. The number of POD  basis obtained using the criterion is mentioned in the bracket of the first column.
From the table, one can see that the eigenvalues corresponding to LHS sampling are better than uniform sampling. In this case, the eigenvalues corresponding to the random sampling and LHS sampling are comparable. In the case of three parameters, the relative error corresponding to uniform sampling is large. 
We report the third eigenvalues in Table~\ref{table:u13_p3}, where we used the first three eigenvectors in the snapshot matrix. Also in these cases, we can see that the results corresponding to the LHS sampling are performing better than the uniform sampling. Moreover, the other two samplings are giving a smaller error than the uniform sampling case.

Hence, for the problem with larger dimensional parameters, sparse sampling seems to perform better than uniform tensor sampling. 
\begin{table}
 	 	\centering
 	\begin{tabular}{|c|c|c|c|c|c|c|c|} 
 		\hline
 	& 
		  {\begin{tabular}[c]{@{}c@{}} $\mu$ \end{tabular}} &
 		 {\begin{tabular}[c]{@{}c@{}} FEM \end{tabular}} & 
 		{\begin{tabular}[c]{@{}c@{}} ROM \end{tabular}}&
 	 	{\begin{tabular}[c]{@{}c@{}}   Error  \end{tabular}}	\\
  \hline
\parbox[t]{2mm}{\multirow{6}{*}{\rotatebox[origin=c]{90}{Random (23)}}} 
&(0.4,0.4,2) &122.8364060    &122.8367075         &$2.4\times 10^{-6}$\\
&(0.4,1.1,2) &70.6551045     &70.6553386         &$3.3\times 10^{-6}$\\
&(1.1,0.4,2) &67.0934953     &67.0936697         &$2.5\times 10^{-6}$\\
&(1.1,1.1,2) &15.0557379     &15.0557756        &$2.5\times 10^{-6}$ \\
&(0.4,0.4,6) &131.7663751    &131.7665266        &$1.1\times 10^{-6}$\\
&(0.4,1.1,6) &79.2294234     &79.2300552       &$7.9\times 10^{-6}$\\
&(1.1,0.4,6) &75.4149311     &75.4159428          &$1.3\times 10^{-5}$\\
&(1.1,1.1,6) &23.1512031     &23.1513903        &$8.0\times 10^{-6}$\\
  \hline
\parbox[t]{2mm}{\multirow{6}{*}{\rotatebox[origin=c]{90}{LHS (21)}}} 
&(0.4,0.4,2) &122.8364060    &122.8365183       &$9.1\times 10^{-7}$\\
&(0.4,1.1,2) &70.6551045     &70.6552179        &$1.6\times 10^{-6}$\\
&(1.1,0.4,2) &67.0934953     &67.0936942        &$2.9\times 10^{-6}$\\
&(1.1,1.1,2) &15.0557379     &15.0558414        &$6.8\times 10^{-6}$ \\
&(0.4,0.4,6) &131.7663751    &131.7665023       &$9.6\times 10^{-7}$\\
&(0.4,1.1,6) &79.2294234     &79.2295204        &$1.2\times 10^{-6}$\\
&(1.1,0.4,6) &75.4149311     &75.4156352        &$9.3\times 10^{-6}$\\
&(1.1,1.1,6) &23.1512031     &23.1513039        &$4.3\times 10^{-6}$\\
  \hline
\parbox[t]{2mm}{\multirow{6}{*}{\rotatebox[origin=c]{90}{Smolyak (19)}}} 
&(0.4,0.4,2) &122.8364060    &122.8376770       &$1.0\times 10^{-5}$\\
&(0.4,1.1,2) &70.6551045     &70.6552022        &$1.3\times 10^{-6}$\\
&(1.1,0.4,2) &67.0934953     &67.1035520        &$1.4\times 10^{-4}$\\
&(1.1,1.1,2) &15.0557379     &15.0577298        &$1.3\times 10^{-4}$ \\
&(0.4,0.4,6) &131.7663751    &131.7674010       &$7.7\times 10^{-6}$\\
&(0.4,1.1,6) &79.2294234     &79.2299386        &$6.5\times 10^{-6}$\\
&(1.1,0.4,6) &75.4149311     &75.4304333        &$2.0\times 10^{-4}$\\
&(1.1,1.1,6) &23.1512031     &23.1556693        &$1.9\times 10^{-4}$\\
  \hline
  \parbox[t]{2mm}{\multirow{6}{*}{\rotatebox[origin=c]{90}{Uniform (20)}}} 
&(0.4,0.4,2) &122.8364060    &122.8376770       &$5.8\times 10^{-5}$\\
&(0.4,1.1,2) &70.6551045     &70.6552022        &$5.1\times 10^{-5}$\\
&(1.1,0.4,2) &67.0934953     &67.1035520        &$3.0\times 10^{-4}$\\
&(1.1,1.1,2) &15.0557379     &15.0577298        &$8.3\times 10^{-4}$ \\
&(0.4,0.4,6) &131.7663751    &131.7674010       &$5.5\times 10^{-5}$\\
&(0.4,1.1,6) &79.2294234     &79.2299386        &$5.1\times 10^{-5}$\\
&(1.1,0.4,6) &75.4149311     &75.4304333        &$6.8\times 10^{-4}$\\
&(1.1,1.1,6) &23.1512031     &23.1556693        &$5.8\times 10^{-4}$\\
  \hline
\end{tabular}
\caption{Comparison of first eigenvalues \eqref{eq:PDE_example1}, obtained using FEM and ROM: the snapshot matrix contains $u_1$ at 27 sample points.}
\label{table:u1_p3}
 \end{table}

\begin{table}
 	 	\centering
 	\begin{tabular}{|c|c|c|c|c|c|c|c|} 
 		\hline
 	& 
		  {\begin{tabular}[c]{@{}c@{}} $\mu$ \end{tabular}} &
 		 {\begin{tabular}[c]{@{}c@{}} FEM \end{tabular}} & 
 		{\begin{tabular}[c]{@{}c@{}} ROM \end{tabular}}&
 	 	{\begin{tabular}[c]{@{}c@{}}   Error  \end{tabular}}	\\
  \hline
\parbox[t]{2mm}{\multirow{6}{*}{\rotatebox[origin=c]{90}{Random (55)}}} 
&(0.4,0.4,2) &331.5127631    &331.5135018       &$2.2\times 10^{-6}$\\
&(0.4,1.1,2) &134.0207799    &134.0267998       &$4.4\times 10^{-5}$\\
&(1.1,0.4,2) &101.6766498    &101.6768535       &$2.0\times 10^{-6}$\\
&(1.1,1.1,2) &39.5356861     &39.5358123        &$3.1\times 10^{-6}$ \\
&(0.4,0.4,6) &341.2134501    &341.2141523       &$2.0\times 10^{-6}$\\
&(0.4,1.1,6) &143.7772626    &143.7862997       &$6.2\times 10^{-5}$\\
&(1.1,0.4,6) &111.3125627    &111.3145985       &$1.8\times 10^{-5}$\\
&(1.1,1.1,6) &48.4323150     &48.4329824        &$1.3\times 10^{-5}$\\
  \hline
 \parbox[t]{2mm}{\multirow{6}{*}{\rotatebox[origin=c]{90}{LHS (51)}}} 
&(0.4,0.4,2) &331.5127631     &331.5130345      &$8.1\times 10^{-7}$\\
 &(0.4,1.1,2) &134.0207799   &134.0208568      &$5.7\times 10^{-7}$\\
  &(1.1,0.4,2) &101.6766498   &101.6775456      &$8.8\times 10^{-6}$\\
  &(1.1,1.1,2) &39.5356861    &39.5369478      &$3.1\times 10^{-5}$ \\
  &(0.4,0.4,6) &341.2134501   &341.2137745       &$9.5\times 10^{-7}$\\
  &(0.4,1.1,6) &143.7772626   &143.7773455      &$5.7\times 10^{-7}$\\
  &(1.1,0.4,6) &111.3125627   &111.3134427      &$7.9\times 10^{-6}$\\
  &(1.1,1.1,6) &48.4323150    &48.4324865      &$3.5\times 10^{-6}$\\
\hline
\parbox[t]{2mm}{\multirow{6}{*}{\rotatebox[origin=c]{90}{Uniform (44)}}} 
  &(0.4,0.4,2) &331.5127631   &331.5167299      &$1.1\times 10^{-5}$\\
  &(0.4,1.1,2) &134.0207799   &134.0493167       &$2.1\times 10^{-4}$\\
  &(1.1,0.4,2) &101.6766498   &101.9762078      &$2.9\times 10^{-3}$\\
  &(1.1,1.1,2) &39.5356861    &39.5726073       &$9.3\times 10^{-4}$ \\
  &(0.4,0.4,6) &341.2134501   &341.2184201     &$1.4\times 10^{-5}$\\
  &(0.4,1.1,6) &143.7772626   &143.8052654      &$1.9\times 10^{-4}$\\
  &(1.1,0.4,6) &111.3125627   &111.6493019      &$3.0\times 10^{-3}$\\
  &(1.1,1.1,6) &48.4323150    &48.4707151      &$7.9\times 10^{-4}$\\
\hline
\parbox[t]{2mm}{\multirow{6}{*}{\rotatebox[origin=c]{90}{Smolyak (45)}}} 
&(0.4,0.4,2) &331.5127631   &331.5163923      &$1.0\times 10^{-5}$\\
 &(0.4,1.1,2) &134.0207799   &134.0217719      &$7.4\times 10^{-6}$\\
  &(1.1,0.4,2) &101.6766498   &101.7104656      &$3.3\times 10^{-4}$\\
  &(1.1,1.1,2) &39.5356861    &39.5563581       &$5.2\times 10^{-4}$ \\
  &(0.4,0.4,6) &341.2134501   &341.2166587      &$9.4\times 10^{-6}$\\
  &(0.4,1.1,6) &143.7772626   &143.7787307      &$1.0\times 10^{-5}$\\
  &(1.1,0.4,6) &111.3125627   &111.35334079      &$3.6\times 10^{-4}$\\
  &(1.1,1.1,6) &48.4323150    &48.4538792      &$4.4\times 10^{-4}$\\
\hline
\end{tabular}
\caption{Comparison of third eigenvalues of \eqref{eq:PDE_example1}, obtained using FEM and ROM: the snapshot matrix contains $u_1,u_2,u_3$ at 27 sample points.}
\label{table:u13_p3}
 \end{table}

\bibliographystyle{plain}
\bibliography{casref}

\end{document}